\documentclass[oneside,11pt]{amsart}

\usepackage{amsmath}
\usepackage{amsthm}

\usepackage{amsfonts}
\usepackage{amssymb}
\usepackage{amsxtra}

\usepackage[arrow, tips, matrix]{xy}
\SelectTips{cm}{12}
\CompileMatrices

\newtheorem{theorem}{Theorem}[section]
\newtheorem{proposition}[theorem]{Proposition}
\newtheorem{lemma}[theorem]{Lemma}
\newtheorem{corollary}[theorem]{Corollary}

\theoremstyle{definition}
\newtheorem{definition}[theorem]{Definition}

\newtheorem{example}[theorem]{Example}

\renewcommand{\bar}[1]{\overline{#1}}
\newcommand{\boundary}{\partial}
\newcommand{\gen}[1]{\langle{#1}\rangle}

\newcommand{\bigpresentation}[2]{ \bigl\langle \, {#1} \bigm| {#2} \,
                                  \bigr\rangle }
\newcommand{\set}[2]{\{\,{#1} \mid {#2} \,\}}
\newcommand{\bigset}[2]{ \bigl\{ \, {#1} \bigm| {#2} \, \bigr\} }

\renewcommand{\setminus}{-}

\newcommand{\field}[1]{\mathbb{#1}}
\newcommand{\Z}{\field{Z}}
\newcommand{\R}{\field{R}}
\newcommand{\N}{\field{N}}
\newcommand{\E}{\field{E}}
\newcommand{\Hyp}{\field{H}}

\newcommand{\inclusion}{\hookrightarrow}
\newcommand{\of}{\circ}
\renewcommand{\hat}{\widehat}

\DeclareMathOperator{\Isom}{Isom}
\DeclareMathOperator{\CAT}{CAT}
\DeclareMathOperator{\Image}{Im}
\DeclareMathOperator{\Aut}{Aut}
\DeclareMathOperator{\Stab}{Stab}
\DeclareMathOperator{\Hull}{Hull}
\newcommand{\ball}[2]{B ( {#1}, {#2} )}
\newcommand{\cball}[2]{\bar{ \ball{#1}{#2} } } 
\newcommand{\nbd}[2]{\mathcal{N}_{#2}({#1})}  
\newcommand{\bignbd}[2]{\mathcal{N}_{#2} \bigl( {#1} \bigr)}
\newcommand{\Set}[1]{\mathcal{#1}}

\usepackage{graphics}
\usepackage{color}
\graphicspath{{pics/}}

\newcommand{\drawrftp}{\begin{center}%
                       \includegraphics{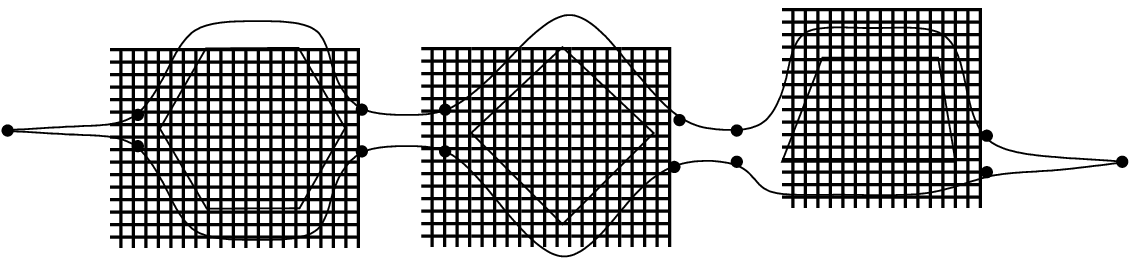}%
                       \end{center}}

\newcommand{\drawLimitSetEFTP}{\begin{center}
                               \begin{picture}(0,0)%
\includegraphics{LimitSetEFTP}%
\end{picture}%
\setlength{\unitlength}{3947sp}%
\begingroup\makeatletter\ifx\SetFigFont\undefined
\def\x#1#2#3#4#5#6#7\relax{\def\x{#1#2#3#4#5#6}}%
\expandafter\x\fmtname xxxxxx\relax \def\y{splain}%
\ifx\x\y   
\gdef\SetFigFont#1#2#3{%
  \ifnum #1<17\tiny\else \ifnum #1<20\small\else
  \ifnum #1<24\normalsize\else \ifnum #1<29\large\else
  \ifnum #1<34\Large\else \ifnum #1<41\LARGE\else
     \huge\fi\fi\fi\fi\fi\fi
  \csname #3\endcsname}%
\else
\gdef\SetFigFont#1#2#3{\begingroup
  \count@#1\relax \ifnum 25<\count@\count@25\fi
  \def\x{\endgroup\@setsize\SetFigFont{#2pt}}%
  \expandafter\x
    \csname \romannumeral\the\count@ pt\expandafter\endcsname
    \csname @\romannumeral\the\count@ pt\endcsname
  \csname #3\endcsname}%
\fi
\fi\endgroup
\begin{picture}(6000,2539)(151,-1844)
\put(2701,-136){\makebox(0,0)[lb]{\smash{\SetFigFont{12}{14.4}{rm}{\color[rgb]{0,0,0}$g$}%
}}}
\put(2776,-1186){\makebox(0,0)[lb]{\smash{\SetFigFont{12}{14.4}{rm}{\color[rgb]{0,0,0}$f$}%
}}}
\put(5326,-1336){\makebox(0,0)[lb]{\smash{\SetFigFont{12}{14.4}{rm}{\color[rgb]{0,0,0}$\eta(k)(x_2)$}%
}}}
\put(3976,539){\makebox(0,0)[lb]{\smash{\SetFigFont{12}{14.4}{rm}{\color[rgb]{0,0,0}$\eta(h)(x_2)$}%
}}}
\put(4276,-1111){\makebox(0,0)[lb]{\smash{\SetFigFont{12}{14.4}{rm}{\color[rgb]{0,0,0}$\alpha_2$}%
}}}
\put(4726,-511){\makebox(0,0)[lb]{\smash{\SetFigFont{12}{14.4}{rm}{\color[rgb]{0,0,0}$\beta_2$}%
}}}
\put(151, 89){\makebox(0,0)[lb]{\smash{\SetFigFont{12}{14.4}{rm}{\color[rgb]{0,0,0}$h(x_1)$}%
}}}
\put(1126,-886){\makebox(0,0)[lb]{\smash{\SetFigFont{12}{14.4}{rm}{\color[rgb]{0,0,0}$\alpha_1$}%
}}}
\put(1726,-661){\makebox(0,0)[lb]{\smash{\SetFigFont{12}{14.4}{rm}{\color[rgb]{0,0,0}$\beta_1$}%
}}}
\put(1576,539){\makebox(0,0)[lb]{\smash{\SetFigFont{12}{14.4}{rm}{\color[rgb]{0,0,0}$g \of \beta_2$}%
}}}
\put(3826,-1636){\makebox(0,0)[lb]{\smash{\SetFigFont{12}{14.4}{rm}{\color[rgb]{0,0,0}$f \of \alpha_1$}%
}}}
\put(1426,-1786){\makebox(0,0)[lb]{\smash{\SetFigFont{12}{14.4}{rm}{\color[rgb]{0,0,0}$k(x_1)$}%
}}}
\put(6151,-1036){\makebox(0,0)[lb]{\smash{\SetFigFont{12}{14.4}{rm}{\color[rgb]{0,0,0}$ $}%
}}}
\end{picture}

                               \end{center}}

\begin{document}

\title{Geometric invariants of spaces
       with isolated flats}

\author{G.~Christopher Hruska$^{\dag}$}
\address{Department of Mathematics\\
  Malott Hall\\
  Cornell University\\
  Ithaca, NY 14853, USA}
\curraddr{Department of Mathematics\\
  University of Chicago\\
  5734 S.~University Ave.\\
  Chicago, IL 60637, USA}
\email{chruska@math.uchicago.edu}
\thanks{$^{\dag}$ Partially supported by a grant from the National Science
  Foundation.}

\begin{abstract}
We study those groups that act properly discontinuously, cocompactly, and
isometrically on $\CAT(0)$ spaces with isolated flats and the
Relative Fellow Traveller Property.  The groups in question include word
hyperbolic $\CAT(0)$ groups as well as geometrically finite Kleinian groups
and numerous $2$--dimensional $\CAT(0)$ groups.
For such a group we show that there is an intrinsic notion of a quasiconvex
subgroup which is equivalent to the inclusion being a quasi-isometric
embedding. We also show that the visual boundary of the $\CAT(0)$ space is
actually an invariant of the group.
More generally, we show that each quasiconvex subgroup of such a group has a
canonical limit set which is independent of the
choice of overgroup.

The main results in this article were established by Gromov and Short in
the word hyperbolic setting and do not extend to arbitrary $\CAT(0)$ groups.
\end{abstract}



\date{\today}

\maketitle

\section{Introduction}
\label{sec:Introduction}

A group is word hyperbolic if it admits a geometric
action (\emph{i.e.}, properly discontinuous, cocompact, and isometric)
on a $\delta$--hyperbolic space.
Numerous geometric features of such an
action have been shown to be
invariants of the group, in particular the visual boundary
(\cite{Gromov87}) and the set of quasiconvex subgroups (\cite{Short91}).
A quasiconvex subgroup of a word hyperbolic group
is again a word hyperbolic group, and its boundary equivariantly
embeds into the boundary of the larger group as a limit set
(\cite{Gromov87}).
These results do not extend from the negatively curved setting to arbitrary
nonpositively curved groups.

The goal of this article is to show that these results do, in fact,
hold for a special class of nonpositively curved groups, namely those
which act geometrically on $\CAT(0)$ spaces with isolated flats
and the Relative Fellow Traveller Property, which were introduced by
the author in \cite{Hruska2ComplexIFP}.
The groups in question include the fundamental group of any compact
nonpositively curved $2$--complex whose $2$--cells are regular Euclidean
hexagons as well as all geometrically finite Kleinian groups.

In the $\delta$--hyperbolic setting Gromov established a Fellow Traveller
Property, which states that quasigeodesics
with common endpoints track close together.  This fundamental property is a
key component of the proofs of the main results in the word
hyperbolic setting.

The Relative Fellow Traveller Property is a generalization of the Fellow
Traveller Property in which pairs of quasigeodesics
fellow travel ``relative to flats'' in a sense that we make precise in
Section~\ref{sec:rftp}.
The Relative Fellow Traveller Property is useful in conjunction with the
isolated flats property,
defined in Section~\ref{sec:IsolatedFlats}, which
roughly states that flat Euclidean subspaces are ``disjoint at infinity.''
Morally speaking, the $\CAT(0)$ spaces with isolated flats
are the $\CAT(0)$ spaces that are closest to being
$\delta$-hyperbolic, while still containing flat subspaces.

Although the notions of isolated flats and the Relative Fellow Traveller
Property were explicitly introduced by the author in
\cite{Hruska2ComplexIFP}, the ideas were implicit in earlier work of
Kapovich--Leeb (\cite{KapovichLeeb95}), Wise
(\cite{Wise96,WiseFigure8}), and Epstein (\cite[Chapter~11]{ECHLPT92}),
and have also been studied by Kleiner.

Let $\rho\colon G \to \Isom (X)$ be a geometric group action.
A subgroup $H$ of~$G$ is \emph{quasiconvex with respect to $\rho$}
if the orbit $Hx$ is a quasiconvex subspace of~$X$ for some basepoint~$x$.
In the word hyperbolic setting
quasiconvexity does not depend on the choice of geometric action~$\rho$
or even on the choice of space~$X$.
In fact, quasiconvexity of~$H$ is equivalent to the intrinsic property
that $H \inclusion G$ is a quasi-isometric embedding.
This result does not extend to the general $\CAT(0)$ setting
(see Section~\ref{sec:Preliminaries}).
Nevertheless, in the presence of isolated flats and the Relative
Fellow Traveller Property, we have the following theorem.

\begin{theorem}[Quasiconvex $\Longleftrightarrow$ undistorted]
\label{thm:QCisQI}
Let $\rho$ be a geometric action of a group~$G$ on a
$\CAT(0)$ space $X$, where $X$ has isolated flats
and the Relative Fellow Traveller Property,
and let $H\le G$ be any finitely generated subgroup.
Then $H$ is quasiconvex with respect to~$\rho$ if and only if
the inclusion $H\inclusion G$ is a quasi-isometric embedding.
\end{theorem}

The boundary of a (complete) $\CAT(0)$ space is the space of all
geodesic rays emanating from a fixed basepoint, endowed
with the compact-open topology.  In the word hyperbolic setting,
the boundary is a group invariant in the sense that,
if a group acts geometrically on two different
$\delta$--hyperbolic spaces then the spaces have the same boundary.
Furthermore, as mentioned above, a quasiconvex subgroup of
a word hyperbolic group is again word hyperbolic, and
there is an equivariant embedding of the boundary of the subgroup
into the boundary of the supergroup as a limit set.

In principle, the same is true in the presence of both isolated flats
and the Relative Fellow Traveller Property.
However, the statement is more subtle
since it is currently unknown whether a quasiconvex subgroup
of a $\CAT(0)$ group is itself $\CAT(0)$.
(Recall that a group is $\CAT(0)$ if it admits a geometric action on a
$\CAT(0)$ space.)

\begin{theorem}[Boundary of a quasiconvex subgroup is
well-defined]
\label{thm:LimitSetWellDefined}
Let $\rho_1$ and~$\rho_2$
be geometric actions of groups $G_1$ and~$G_2$ on $\CAT(0)$ spaces
$X_1$ and~$X_2$ each having isolated flats and
the Relative Fellow Traveller Property.
For each $i$, let $H_i\le G_i$ be a quasiconvex subgroup with respect
to~$\rho_i$.
Then any isomorphism $\eta \colon H_1 \to H_2$
induces an $\eta$--equivariant homeomorphism
$\Lambda H_1 \to \Lambda H_2$.
\end{theorem}

Roughly speaking the idea is that, if a group $H$ is a quasiconvex
subgroup of two different groups $G_1$ and~$G_2$,
then the limit set of~$H$ is the same in the boundary of both groups.
If $H$ is itself a $\CAT(0)$ group, then this limit set must also be the
boundary of~$H$.

An immediate corollary is that the boundary is a group invariant
for $\CAT(0)$ groups with isolated flats and the Relative Fellow Traveller
Property.  Croke and Kleiner showed that this corollary does not extend to
the general $\CAT(0)$ setting (\cite{CrokeKleiner00}).
In fact Wilson has shown that the Croke--Kleiner construction
produces a continuous family of homeomorphic $2$--complexes
whose universal covers all have topologically distinct boundaries
(\cite{WilsonBoundary}).

In order to prove the main theorems, we first prove some basic
algebraic facts due to Wise (personal communication)
about geometric actions on $\CAT(0)$ spaces with
isolated flats in Section~\ref{sec:IsolatedFlats}.
In particular, we show that maximal flats correspond to maximal
virtually abelian subgroups of rank at least two
in Theorem~\ref{thm:PeriodicFlats}.  We also establish that these
virtually abelian subgroups lie in only finitely many conjugacy classes
(Theorem~\ref{thm:FinConjClassesAbelians}).
One notable corollary to this analysis is the following.

\begin{theorem}[$\Z \times \Z$ subgroups]
Suppose $G$ acts geometrically on a $\CAT(0)$ space with isolated flats.
Then either $G$ is word hyperbolic or $G$ contains a $\Z\times \Z$ subgroup.
\end{theorem}

We prove Theorem~\ref{thm:QCisQI} in
Section~\ref{sec:QCWellDefined}.  The proof is a more detailed version
of the techniques used to prove Theorem~\ref{thm:PeriodicFlats}.

Then in Section~\ref{sec:Quasiflats} we examine the geometry of \emph{quasiflats},
\emph{i.e.}, the quasi-isometrically embedded Euclidean subspaces
of dimension at least two.  In particular
we prove the following theorem, which generalizes a lemma proved by
Schwartz in his study of quasi-isometric rigidity of nonuniform lattices
in rank one symmetric spaces (\cite{Schwartz95}).

\begin{theorem}[Quasiflats are close to flats]\label{thm:Quasiflat}
Let $X$ be a $\CAT(0)$ space with isolated flats and the
Relative Fellow Traveller Property.
Given constants $\lambda$ and~$\epsilon$, there is a constant $D =
D(\lambda,\epsilon,X)$ such that each $(\lambda,\epsilon)$--quasiflat
lies in a $D$--neighborhood of some flat~$F$.
\end{theorem}

Theorem~\ref{thm:Quasiflat} and the Flat Torus Theorem are key components in
the proof of the following theorem,
which improves the Relative Fellow Traveller Property to a genuine Fellow
Traveller Property under an equivariance assumption.

\begin{theorem}[Equivariant Fellow Traveller Property]
\label{thm:EquivariantFTP}
Let $G$ act geometrically on two $\CAT(0)$ spaces $X$ and~$Y$.
Suppose further that $X$ has isolated flats and the Relative
Fellow Traveller Property.
Then any $G$--equivariant quasi-isometry $X \to Y$ maps geodesics in~$X$
uniformly close to geodesics in~$Y$.
\end{theorem}

As an immediate corollary, we get an alternate proof of the group
invariance of quasiconvexity without using the techniques
of Section~\ref{sec:QCWellDefined}.
We also obtain a direct proof of the invariance of the boundary.
To prove the more general result of Theorem~\ref{thm:LimitSetWellDefined}
requires a more detailed analysis which combines techniques
from the proofs of Theorems \ref{thm:QCisQI}
and~\ref{thm:EquivariantFTP}.  We undertake this analysis in
Section~\ref{sec:LimitSets}.

In \cite{Hruska2ComplexIFP}, the author shows that among proper,
compact $\CAT(0)$ $2$--complexes the isolated flats property is equivalent
to the Relative Fellow Traveller Property.  Furthermore, Wise has shown that
these equivalent properties are satisfied if and only if the $2$--complex
does not contain an isometrically embedded triplane
(see \cite{Hruska2ComplexIFP} for a proof).
A triplane is the space obtained by gluing three Euclidean halfplanes
together along their boundary line.
For instance the universal cover of a compact nonpositively curved
$2$--complex whose $2$--cells are all regular Euclidean hexagons
cannot contain a triplane, and hence has isolated flats and the
Relative Fellow Traveller Property.

Every geometrically finite
Kleinian group acts geometrically on a truncated version of the convex
hull of the limit set.  It is a well-known folk theorem that this
``truncated convex hull'' is a $\CAT(0)$ space.
In Section~\ref{sec:GeomFinite} we prove the following result.

\begin{theorem}
\label{thm:GeomFiniteRFTP}
The truncated convex hull associated to a geometrically finite
subgroup $\Gamma \le \Isom(\Hyp^n)$ has isolated flats and the
Relative Fellow Traveller Property.
\end{theorem}

The Relative Fellow Traveller Property is established using a
technical result of Epstein (\cite[11.3.1]{ECHLPT92}).

Based on the evidence above,
it seems likely that the isolated flats property and the Relative
Fellow Traveller Property are equivalent for arbitrary
(proper and cocompact) $\CAT(0)$ spaces.

The author has been informed that Kleiner has unpublished work from
1997 related to the paper \cite{CrokeKleiner02}
which also proves some of the results in this article.
In particular he showed that equivariant quasi-isometries between
CAT(0) spaces with isolated flats map geodesics to within
uniform Hausdorff distance of geodesics, which implies that
the spaces have a well-defined boundary.

\subsection{Acknowledgements}
\label{subsec:Acknowledgements}

The results in this article were originally published as part of my Ph.D.
dissertation at Cornell University.  That dissertation was prepared under the
guidance of Daniel Wise and Karen Vogtmann.
I would like to thank Karen for supporting my desire to work with Dani on
this project, although he was a postdoc at the time.

Dani originally conjectured several of the main results of this article
and had established the results of Section~\ref{sec:IsolatedFlats} in the two
dimensional setting before I started this project.  In fact a variant of
Theorem~\ref{thm:PeriodicFlats} appeared as Proposition~4.0.4 in Dani's own
thesis \cite{Wise96}.
I thank him for his vision and encouragement.

During this project,
I benefited from numerous conversations with Marshall Cohen, Jon
McCammond, John Meier, Kim Ruane, and surely others that I have forgotten to
mention.  I would also like to thank the referee for many useful comments
that I hope have improved the exposition of this paper.

\section{Geometric preliminaries}
\label{sec:Preliminaries}

This section is a review of some basic geometric facts that we will need
throughout this article.  It also serves to establish notation.
A good reference for the facts discussed here is
\cite{BH99}.

We use two distinct metrics on the set of subsets
of a metric space~$X$.
If $A$ and~$B$ are subsets of~$X$, the \emph{distance}
between $A$ and~$B$ is defined by
\[
   d(A,B) = \inf \,\bigset{d(a,b)}{a\in A, \ b\in B}.
\]
The \emph{Hausdorff distance} between $A$ and~$B$ is
\[
  d_H(A,B) = \inf \, \bigset{\epsilon}{A \subseteq \nbd{B}{\epsilon}
   \text{ and } B \subseteq \nbd{A}{\epsilon}},
\]
where $\nbd{C}{\epsilon}$ denotes the $\epsilon$--neighborhood of~$C$.

\subsection{Quasi-isometries}
\label{subsec:QuasiIsometries}

\begin{definition}[Quasi-isometry]
Let $X$ and~$Y$ be metric spaces.
A \emph{$(\lambda,\epsilon)$--quasi-isometric embedding}
of~$X$ into~$Y$ is a function
$f \colon X \to Y$ satisfying
\[
   \frac{1}{\lambda}\, d(a,b) - \epsilon  \le  d \bigl( f(a),f(b) \bigr)
     \le \lambda \, d(a,b) + \epsilon
\]
for all $a,b \in X$.
If, in addition, every point of~$Y$ lies in the $\epsilon$--neighborhood
of the image of~$f$, then $f$ is a
\emph{$(\lambda,\epsilon)$--quasi-isometry}
and $X$ and~$Y$ are
\emph{quasi-isometric}.
\end{definition}

Every quasi-isometry $f$ has a \emph{quasi-inverse} $g$
with the property that the maps $f\of g$ and $g\of f$
are each within a bounded distance of the identity.

A \emph{geodesic} in a metric space~$X$ is an isometric embedding
$I \to X$, where $I \subseteq \R$ is an interval.
A metric space is \emph{geodesic} if every pair of points is connected
by a geodesic.
A group action is \emph{geometric} if it is
properly discontinuous, cocompact, and isometric.
The following well-known result was discovered by
Efromovich and \v{S}varc (\cite{Efromovich53,Svarc55})
and rediscovered by Milnor (\cite{Milnor68}).

\begin{theorem}\label{thm:SvarcMilnor}
Suppose a group~$G$ acts geometrically on a geodesic space~$X$.
Then the map $G \to X$ given by $g \mapsto g(x_0)$ for some $x_0\in X$
is a quasi-isometry.
\end{theorem}

\subsection{CAT(0) spaces}
\label{subsec:CAT0}

\begin{definition}
Given a geodesic triangle~$\Delta$ in~$X$, a \emph{comparison triangle}
for $\Delta$ is a triangle in the Euclidean plane with the
same edge lengths as~$\Delta$.  A geodesic space is $\CAT(0)$
if distances between points on any geodesic triangle~$\Delta$
are less than or
equal to the distances between the corresponding points on a
comparison triangle for~$\Delta$.
\end{definition}

\begin{theorem}[Convexity of the $\CAT(0)$ metric]
Let $\gamma$ and $\gamma'$ be geodesic segments in a $\CAT(0)$ space,
each parametrized from $0$ to~$1$ proportional to arclength.
Then for each $t\in [0,1]$ we have
\[
   d\bigl( \gamma(t),\gamma'(t) \bigr)
   \le (1-t) \, d\bigl( \gamma(0),\gamma'(0) \bigr)
      + t \, d\bigl( \gamma(1),\gamma'(1) \bigr).
\]
\end{theorem}

\begin{theorem}[Orthogonal projection]
Let $C$ be a complete, convex subspace of a $\CAT(0)$ space~$X$.
Then there exists a unique map $\pi \colon X \to C$,
called the \emph{orthogonal projection} of~$X$ onto~$C$,
satisfying the following properties.
\begin{enumerate}
   \item For each $x \in X$, we have $d \bigl( x, \pi(x) \bigr) = d(x,C)$.
   \item For every $x$, $y \in X$ we have
      $d(x,y) \ge d \bigl( \pi(x), \pi(y) \bigr)$.
\end{enumerate}
\end{theorem}

Recall that an isometry~$g$ of a metric space~$X$ is \emph{semisimple} if
some point of~$X$ is moved a minimal distance by~$g$.  If $G$ acts
geometrically on a metric space~$X$, then every element of~$G$
is a semisimple isometry of~$X$ (see \cite[II.6.10(2)]{BH99}).

\begin{theorem}[Flat Torus Theorem]\label{thm:FlatTorus}
Let $A$ be a free abelian group of rank~$k$ acting properly
discontinuously by semisimple
isometries on a $\CAT(0)$ space~$X$.  Then $A$ stabilizes some
$k$--flat~$F$,
and the action of~$A$ on~$F$ is by Euclidean translations with quotient a
$k$--torus.
\end{theorem}

\subsection{Quasiconvexity}
\label{subsec:Quasiconvexity}

\begin{definition}[Quasiconvex subspace]\label{def:QCsubspace}
A subspace~$Y$ of a geodesic metric space~$X$
is \emph{$\nu$--quasi\-convex}, for $\nu \ge 0$, if every
geodesic in~$X$ connecting two points of~$Y$ lies inside a
$\nu$--neighborhood of~$Y$.
The subspace~$Y$ is \emph{quasiconvex} if there exists a nonnegative
constant~$\nu$ such that $Y$~is $\nu$--quasiconvex.
\end{definition}

\begin{definition}[Quasiconvex subgroup]
Let $\rho\colon G\to \Isom{X}$ be a geometric action of the group~$G$ on the
$\CAT(0)$ space~$X$.
A subgroup~$H$ of~$G$ is \emph{quasiconvex}
with respect to~$\rho$ if there is a point $x_0\in X$ such that the
orbit $Hx_0$ is a quasiconvex subspace of~$X$.
\end{definition}

\begin{theorem}\label{thm:AlonsoBridson}
Quasiconvex subgroups are finitely generated and
quasi-isometrically embedded.
Furthermore, if two subgroups are quasiconvex with respect to
the same action, then their intersection is again
quasiconvex.
\end{theorem}

In the general $\CAT(0)$ setting,
quasiconvexity depends on the choice of action,
as illustrated in the following example.

\begin{example}\label{exmp:QCExmp}
Consider the group
\[
   G = F_2 \times \Z = \gen{a,b} \times \gen{t}
     = \bigpresentation{a,b,t}{[a,t], [b,t]},
\]
and let $\rho \colon G \to \Isom(X)$ be the natural action
of~$G$ on the universal cover of the presentation $2$--complex,
metrized as a product of two trees.
Since $\rho$ respects this product decomposition,
the direct factor $H = \gen{a,b}$ is quasiconvex
with respect to~$\rho$.

Observe that $H' = \gen{a,bt}$ is not quasiconvex
with respect to~$\rho$.
For if it were, then $H \cap H'$ would be finitely generated.
But $H\cap H'$ is the subgroup of all elements in the free group~$H$
for which the exponent sum of~$b$ is zero, which is not finitely
generated.

The map $\phi \in \Aut(G)$ given by
\[
  a \mapsto a \qquad b \mapsto bt \qquad t \mapsto t
\]
sends $H$ to~$H'$.
So the $\rho$~action of~$H'$ is the same as the $\rho\of\phi$
action of~$H$.  In other words, $H$ is not quasiconvex with respect to
$\rho\of\phi$.
\end{example}

\section{Isolated flats}
\label{sec:IsolatedFlats}

In this section, we define the notion of a $\CAT(0)$ space with
\emph{isolated flats}.
We prove some
algebraic facts about groups which act geometrically on
$\CAT(0)$ spaces with isolated flats.
In particular, we show
in Theorem~\ref{thm:PeriodicFlats} that
maximal flats are in one-to-one correspondence with maximal free abelian
subgroups of rank at least two.
We also prove
Theorem~\ref{thm:FinConjClassesAbelians}, which states that
such groups have finitely many conjugacy classes of
maximal virtually abelian subgroups of rank at least two.

The results in this section were proved by Wise in the
$2$--dimensional setting (personal communication).
The proofs here are straightforward generalizations of those given by Wise.
In fact, a variant of Theorem~\ref{thm:PeriodicFlats} appears as
Proposition~4.0.4 in \cite{Wise96}.

\begin{definition}[Flats]\label{def:Flats}
A \emph{flat} in a $\CAT(0)$ space~$X$ is an isometric embedding
of Euclidean space~$\E^k$ into~$X$ for some $k \ge 2$.
A $k$--flat is a flat of dimension~$k$.
\end{definition}

\begin{definition}[Isolated flats]\label{def:IsolatedFlats}
A $\CAT(0)$ space~$X$ has \emph{isolated flats}
if it contains a family~$\Set{F}$ of flats
with the following properties.
\begin{enumerate}
\item\label{item:IFPMaximal}
   (Maximal) There is a constant~$B$ such that every flat~$F$ in~$X$
   is contained in a $B$--neighborhood of some
   flat $F' \in \Set{F}$.
\item\label{item:IFPIsolated}
   (Isolated) There is a function $\psi \colon \R_+ \to \R_+$
   such that for every pair of distinct flats $F_1, F_2 \in \Set{F}$
   and for every $k\ge 0$, the intersection $\nbd{F_1}{k} \cap \nbd{F_2}{k}$
   of $k$--neighborhoods of $F_1$ and~$F_2$ has diameter at most $\psi(k)$.
\item\label{item:IFPEquivariant}
   (Equivariant) The set of flats~$\Set{F}$ is invariant under the action
   of $\Isom(X)$.
\end{enumerate}
\end{definition}

Observe that $\delta$--hyperbolic $\CAT(0)$ spaces vacuously satisfy
the isolated flats property since such spaces do not contain flats.

We note the following immediate consequence of isolated flats,
which will be useful in the sequel.

\begin{lemma}\label{lem:IFPEquivalent}
Let $X$ be a $\CAT(0)$ space with isolated flats.
For every $k \ge 0$, each flat disc~$D$ in~$X$ of radius
at least $\psi(k)$
lies in a $k$--neighborhood of at most one flat $F \in \Set{F}$.
\qed
\end{lemma}

The following proposition shows that in any proper $\CAT(0)$ space with
isolated flats, the family~$\Set{F}$
is locally finite.

\begin{proposition}[Locally finite]\label{prop:LocallyFinite}
Let $X$ be a proper $\CAT(0)$ space with isolated flats.
Then only finitely many flats from the family~$\Set{F}$
intersect any compact set $K \subseteq X$.
\end{proposition}

\begin{proof}
It suffices to show that only finitely many flats in~$\Set{F}$
intersect any metric ball $\ball{x_0}{r}$.
By Lemma~\ref{lem:IFPEquivalent}, we can choose a constant~$R$
sufficiently large that every flat disc~$D$ in~$X$ of radius at least~$R$
lies in a $1$--neighborhood of at most one flat $F \in \Set{F}$.
Let $\{F_i\}$ be the collection of all flats $F_i \in \Set{F}$ which
intersect
the ball $\ball{x_0}{r}$.  Let $p_i$ be the point in $F_i$ closest to~$x_0$,
and let $D_i$ be the closed disc of radius~$R$ in~$F_i$ centered at~$p_i$.
Note that every such disc lies inside the closed ball $\cball{x_0}{r+R}$,
which is compact since $X$ is proper.

Suppose by way of contradiction that the collection $\{F_i\}$ is infinite.
Passing to a subsequence if necessary,
we may assume that the sequence of discs $\{D_i\}$ converges
in the Hausdorff metric
(see, for instance, \cite[I.5.31]{BH99}).
In particular, it is a Cauchy sequence with respect to the metric $d_H$,
so some pair
$D_i,D_j$ with $i\ne j$ has $d_H(D_i,D_j) < 1$.  But then $D_i$ lies inside
a $1$--neighborhood of the distinct flats $F_i$ and $F_j$,
contradicting our choice of~$R$.
\end{proof}

\begin{corollary}\label{cor:FinConjClassesFlats}
Suppose a group~$G$ acts geometrically on a $\CAT(0)$ space~$X$
with isolated flats.  Then
$G$ contains only finitely many conjugacy classes of stabilizers of
flats $F \in \Set{F}$.
\end{corollary}

\begin{proof}
Fix a compact set $K$ in~$X$ whose $G$--translates cover $X$.
Every flat $F \in \Set{F}$ intersects $g(K)$ for some $g \in G$.
So the flat $g^{-1}(F)$ intersects~$K$ and has a stabilizer conjugate to
the stabilizer of~$F$.  Since only finitely many flats in~$\Set{F}$
intersect~$K$, we see that there are only finitely many conjugacy classes of
flat stabilizers.
\end{proof}

\begin{definition}[Periodic]\label{def:periodic}
Suppose a group $G$ acts geometrically on a metric space $X$.
A $k$--flat~$F$ in~$X$ is \emph{periodic}
if there is a free abelian subgroup $A \le G$
of rank~$k$ that acts by translations on~$F$ with quotient a $k$--torus.
\end{definition}

The following theorem is, in a sense, a converse to the Flat Plane Theorem
in the context of isolated flats.

\begin{theorem}[Flats are periodic]\label{thm:PeriodicFlats}
Suppose a group~$G$ acts geometrically on a $\CAT(0)$ space~$X$
which has isolated flats.
Then every $F \in \Set{F}$ is periodic.
\end{theorem}

\begin{proof}
Since $G$ acts cocompactly, the quotient $G \backslash X$
has a bounded diameter~$r$.
Choose a $k$--flat $F \in \Set{F}$, and
let $\set{g_j}{j \in \N}$ be a minimal set of group elements
such that every point of~$F$ lies within a distance~$r$ of some
$g_j(x)$, where $x$ is a fixed basepoint in~$X$.
Then the flat $F_j = g_j^{-1}(F)$ intersects
$\cball{x}{r}$.
Since $\Set{F}$ is invariant under isometries of~$X$,
each flat~$F_j$ is an element of~$\Set{F}$.
But only finitely many flats in~$\Set{F}$ intersect this ball.
So the collection $\set{F_j}{j \in \N}$ is finite.

Let $G_F$ denote the stabilizer of~$F$.  Notice that if two flats $F_i$
and $F_j$
coincide, then $g_j g_i^{-1}$ is an element of~$G_F$.
For each $j$, let $x_j$ denote the point in $F'$ closest to $g_j(x)$.
Then the $x_j$ lie in only finitely many different
$G_F$--orbits.  So for some~$j$ the $G_F$--orbit of $x_j$ is infinite
and does not lie inside a bounded neighborhood of any hyperplane of~$F'$.
Since $G_F$ acts properly discontinuously by isometries on the Euclidean
space~$F$, it follows that $G_F$ has a free abelian subgroup~$A$ of finite
index and finite rank (see Corollary~4.1.13 of \cite{Thurston97}).
Then by Theorem~\ref{thm:FlatTorus}
there is an $m$--flat $F_A$ in~$F$ stabilized by~$G_F$
on which $A$ acts by translations,
where $m$ is equal to the rank of~$A$.
So each orbit under~$G_F$
lies within a bounded distance of~$F_A$.  (Although the specific bound
depends
on the choice of orbit.)
Since $F$ contains a $G_F$--orbit which does not lie in a bounded
neighborhood
of any hyperplane, we must have $m=k$, in other words, $F_A = F$.
Since $A$ acts freely and cocompactly on~$F$ by translations,
the quotient $A \backslash F$ is a $k$--torus as desired.
\end{proof}

The following algebraic result is an immediate consequence of
Theorem~\ref{thm:PeriodicFlats}.

\begin{corollary}\label{cor:Flat=Group}
Suppose a group~$G$ acts geometrically on a $\CAT(0)$ space~$X$
with isolated flats.
Then $X$ contains
a $k$--flat if and only if\/ $G$ contains a subgroup isomorphic to~$\Z^k$.
\qed
\end{corollary}

Taken together, Corollary~\ref{cor:FinConjClassesFlats} and
Theorem~\ref{thm:PeriodicFlats} have the following algebraic consequence.

\begin{theorem}\label{thm:FinConjClassesAbelians}
If $G$ acts geometrically on a $\CAT(0)$ space~$X$ with isolated flats,
then $G$ contains only finitely many conjugacy classes of
maximal virtually abelian subgroups of rank at least two.
\end{theorem}

\begin{proof}
By Corollary~\ref{cor:FinConjClassesFlats}, it suffices to
show that the set of all stabilizers of flats in~$\Set{F}$ is the same
as the set~$\Set{A}$ of
all maximal virtually abelian subgroups of rank at least two of~$G$.

By the Flat Torus Theorem,
each $A\in \Set{A}$ stabilizes a flat~$E$.
But $E$ lies in a tubular neighborhood of a
\emph{unique} flat $F \in \Set{F}$.
The equivariance of~$\Set{F}$ shows that
$A$ is contained in the virtually abelian group stabilizing~$F$,
so in fact $A = \Stab(F)$.
Conversely, for each $F \in \Set{F}$,
Theorem~\ref{thm:PeriodicFlats} gives
$\Stab(F) \subseteq A = \Stab(F')$ for some $A \in \Set{A}$ and
$F' \in \Set{F}$.  By isolated flats, we must have $F=F'$.
\end{proof}

\section{Invariance of quasiconvexity}
\label{sec:QCWellDefined}

A key step in
the proof of Theorem~\ref{thm:QCisQI} is a generalization
of the proof of Theorem~\ref{thm:PeriodicFlats}.
In the earlier proof, we considered a flat coarsely covered
by orbit points $g_i(x_0)$.  We used elements of the form
$g_j g_i^{-1}$ to generate a large virtually abelian
group stabilizing the given flat.

In this section, we prove the following lemma
involving a curve in a flat, which is coarsely covered by orbit
points $h_i(x_0)$ in a subgroup~$H$.  We will use elements of the form
$h_j h_i^{-1}$ to generate a virtually abelian
subgroup of~$H$ that densely fills a subflat coarsely containing
the given curve.

\begin{lemma}\label{lem:PushAcrossOneFlat}
Suppose a group~$G$ acts geometrically on a $\CAT(0)$ space~$X$
with isolated flats.
For each constant $\mu > 0$,
there is a positive constant $L=L(\mu)$ having the following property.
Let $\alpha\colon [0,1]\to F$ be any path in a flat $F \in \Set{F}$,
and let $H$ be a subgroup of~$G$.
Suppose $\Image(\alpha)$ lies in a $\mu$--neighborhood of some orbit
$Hx$ under~$H$.
Then:
\begin{enumerate}
  \item $\Image(\alpha)$ lies in an $L$--neighborhood of a flat
        subspace $\hat{F}$ of\/~$F$ on which a free abelian subgroup
        $\hat{B} \le H$ acts cocompactly by translations.
  \item For each $y \in \hat{F}$ the orbit $\hat{B}y$ is $L$--dense
        in~$\hat{F}$.
  \item The geodesic segment~$\gamma$ connecting $\alpha(0)$ and $\alpha(1)$
        lies in an $L$--neighborhood of\/~$Hx$.
\end{enumerate}
\end{lemma}

The proof of Lemma~\ref{lem:PushAcrossOneFlat} uses the following
elementary lemma, whose proof we leave as an exercise.

\begin{lemma}\label{lem:Orbits}
Suppose a group~$G$ acts by isometries on a metric space~$X$.
If some connected set~$C$ in~$X$ lies in a $\kappa$--neighborhood of the
union of $n$ distinct $G$--orbits, then $C$ lies in a
$(2\kappa n)$--neighborhood of a single orbit. \qed
\end{lemma}

\begin{proof}[Proof of Lemma~\ref{lem:PushAcrossOneFlat}]
Given a path~$\alpha \colon [0,1] \to F$ in some flat $F\in\Set{F}$,
choose a minimal set $\set{h_i}{i \in I}$ of elements of~$H$
so that every point of $\Image(\alpha)$ lies within a distance~$\mu$ of some
$h_i(x)$.
For convenience, replace the points $h_i(x)$ with
points inside~$F$ as follows.  For each $i \in I$,
let $x_i = \pi \bigl( h_i(x) \bigr)$, where $\pi \colon X \to F$ is
the orthogonal projection onto~$F$.
Since projections do not increase distances, every point of $\Image(\alpha)$
lies within a distance~$\mu$ of some~$x_i$.

For each $i \in I$, the flat $h_i^{-1}(F)$ is an element of~$\Set{F}$
intersecting the ball $\ball{x}{\mu}$,
and $h_i^{-1}(x_i)$ is the closest point in this flat to the basepoint~$x$.
Let $N = N(\mu)$ be the number of flats in~$\Set{F}$ which intersect
this ball, which we know to be finite by Proposition~\ref{prop:LocallyFinite}.
Then the set $\bigset{h_i^{-1}(x_i)}{i \in I}$ contains at most $N$ elements.
If two flats $h_i^{-1}(F)$ and $h_j^{-1}(F)$ coincide,
then $h_j h_i^{-1}$ stabilizes~$F$ and maps $x_i$ to~$x_j$.
If we let $H_F$ denote the elements of~$H$ which stabilize~$F$,
then the points~$x_i$ lie in at most~$N$ distinct $H_F$--orbits.

By Theorem~\ref{thm:PeriodicFlats}, we know that
$G_F$, the subgroup of~$G$ stabilizing~$F$, has a free abelian finite index
subgroup~$A_F$ which acts on~$F$ by Euclidean translations.
But there are only finitely many conjugacy classes of the stabilizers~$G_F$
for $F \in \Set{F}$
by Corollary~\ref{cor:FinConjClassesFlats}.
So there is a universal bound~$M$ on the index $[G_F : A_F]$.
Consequently, $H_F$ also has a free abelian subgroup~$B$
of index at most~$M$ which acts by translations on~$F$.
So the points $\{x_i\}$ lie in at most $MN$ distinct $B$--orbits.
Since $\Image(\alpha)$ is connected, Lemma~\ref{lem:Orbits}
shows that $\Image(\alpha)$
lies in a $\mu'$--neighborhood of a single $B$--orbit,
where $\mu' = 2\mu MN$.

Choose a point $y$ in~$F$ and a collection $\set{b_j}{j \in J}$
in~$B$ so that
\[
   d \bigl( b_j(y), b_{j+1}(y) \bigr) < 2\mu'
\]
and $\Image(\alpha)$ lies in a $\mu'$--neighborhood of
$\bigcup_j \bigl\{ b_j(y) \bigr\}$.
Let $\hat{B}$ be the subgroup of~$B$ generated by all elements of the form
$b_{j+1} b_j^{-1}$ for $j \in J$.
Then $\hat{B}$ stabilizes some $k$--flat $\hat{F}$ containing~$y$,
where $k$ is the rank of~$\hat{B}$.

Since the orbit $\hat{B}y$ lies entirely within~$\hat{F}$,
it follows that $\Image(\alpha)$ lies inside a $\mu'$--neighborhood
of~$\hat{F}$.
Furthermore, since $\hat{B}$ is generated by elements with translation
length at most $2\mu'$,
the $k$--flat $\hat{F}$ lies in a $2k\mu'$--neighborhood of $\hat{B}y$.
Since the rank~$k$ of~$\hat{B}$ is bounded by the dimension of the
largest flat in~$X$, we see that $\Image(\alpha)$ lies within a uniformly
bounded
neighborhood of the orbit $\hat{B}y$.

As the $\mu'$--neighborhood of~$\hat{F}$ is convex, the geodesic~$\gamma$
connecting the endpoints of~$\alpha$ also lies
uniformly close to~$\hat{B}y$,
so $\Image(\gamma)$ lies within a uniform neighborhood
of the full orbit $Hy$.
Since each endpoint of~$\gamma$ also lies within a distance~$\mu$
of the original
orbit $Hx$, we see that the Hausdorff distance between the orbits
$Hx$ and $Hy$ is uniformly bounded.  So $\Image(\gamma)$ lies in
a uniformly bounded neighborhood of~$Hx$ as desired.
\end{proof}

\begin{lemma}\label{lem:PushAcrossFlats}
Let $X$ be a $\CAT(0)$ space
with isolated flats
and the Relative Fellow Traveller Property.
Fix constants $\mu$,~$\lambda$, and~$\epsilon$.  Then there exists
a positive constant~$\nu$ so that the following property holds.

Suppose a group~$G$ acts geometrically on~$X$.
Let $\alpha\colon [0,1]\to X$
be a $(\lambda, \epsilon)$--quasigeodesic, and let $\gamma$ be the geodesic
connecting $\alpha(0)$~and~$\alpha(1)$.
Let $H\le G$ be a subgroup
such that $\Image(\alpha)$ lies in a $\mu$--neighborhood of\/~$Hx$
for some basepoint~$x\in X$.
Then $\Image(\gamma)$ lies in a $\nu$--neighborhood of\/~$Hx$.
\end{lemma}

\begin{proof}
By the Relative Fellow Traveller Property
we know that $\alpha$ and~$\gamma$ track $\delta$--close
relative to some sequence of flats in~$\Set{F}$.
The result is clear for any subsegment of~$\gamma$ which lies within a
$\delta$--neighborhood of $\Image(\alpha)$.
So we only need to verify the result for pieces of~$\gamma$
which wander away from $\alpha$.

Let $\xi$ be a subpath of~$\alpha$ whose endpoints are within
a distance~$\delta$
of $\Image(\gamma)$
and which stays in a $\delta$--neighborhood of some flat $F \in \Set{F}$.
Letting $\pi\colon X \to F$ denote the orthogonal projection,
the Hausdorff distance between $\Image(\xi)$ and $\Image(\pi\of\xi)$ is
at most~$\delta$,
so $\Image(\pi\of\xi)$ lies in a $(\delta+\mu)$--neighborhood of $Hx$.

Applying Lemma~\ref{lem:PushAcrossOneFlat} to the curve $\pi\of\xi$,
we get a constant $L$ so that the geodesic~$\eta$ connecting
the endpoints of
$\pi\of\xi$ lies in an $L$--neighborhood of $Hx$.
The result follows from the observation that
the endpoints of~$\eta$ are within a distance $2\delta$ of $\Image(\gamma)$.
\end{proof}

At this point, the proof of Theorem~\ref{thm:QCisQI} is nearly immediate.

\begin{proof}[Proof of Theorem~\ref{thm:QCisQI}]
The direction $(\Rightarrow)$ is Theorem~\ref{thm:AlonsoBridson}.
For $(\Leftarrow)$, suppose $H\inclusion G$ is a quasi-isometric embedding.
Then there is an $H$--equivariant $(\lambda,\epsilon)$--quasi-isometric
embedding
$\phi \colon Y \to X$ where $Y$ is some Cayley graph for~$H$.
Let $x=\phi(1)$.

Note that $H$ is a $(1/2)$--quasiconvex subspace of~$Y$.
Let $\gamma$ be the geodesic in~$Y$ joining two arbitrary points
$h_1$~and~$h_2$ of~$H$.
Then $\phi\of \gamma$ is a $(\lambda,\epsilon)$--quasigeodesic
joining $h_1(x)$ and~$h_2(x)$ and lying in a
$({\lambda}/{2} + \epsilon)$--neighborhood of~$Hx$.
It now follows from Lemma~\ref{lem:PushAcrossFlats}
that there is a constant~$\nu=\nu(\lambda,\epsilon)$
so that the geodesic $\bigl[ h_1(x),h_2(x) \bigr]$ lies in a
$\nu$--neighborhood of~$Hx$.  Thus $H$~is a $\nu$--quasiconvex subgroup.
\end{proof}

\section{The Relative Fellow Traveller Property}
\label{sec:rftp}

\begin{figure}
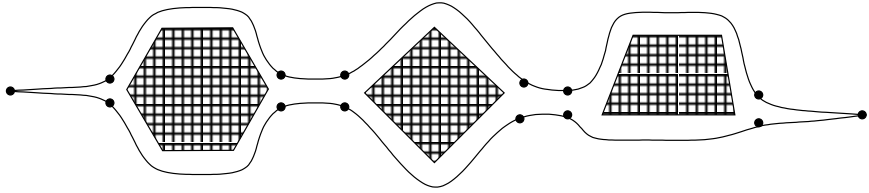

\drawrftp
\caption{A pair of paths which fellow travel relative to flats}
\label{fig:RelativeFTP}
\end{figure}

Roughly speaking, the idea is that the two paths alternate between tracking
close together and travelling near a common flat as illustrated in
Figure~\ref{fig:RelativeFTP}.

\begin{definition}[Fellow travelling relative to flats]
\label{def:RelativeFellowTravelling}
A pair of paths
\[
   \alpha\colon [0,a]\to X \quad \text{and}
   \quad \alpha'\colon [0,a']\to X
\]
in
a space \emph{$L$--fellow travel relative to a sequence of flats}
$(F_1,\dots,F_n)$
if there are partitions
\[
  0 = t_0 \le s_0 \le t_1 \le s_1 \le \dots \le t_n \le s_n = a
\]
and
\[
  0 = t'_0\le s'_0\le t'_1\le s'_1\le \dots \le t'_n\le s'_n= a'
\]
so that for $0\le i \le n$ the Hausdorff distance between the
sets $\alpha \bigl( [t_i,s_i] \bigr)$ and
$\alpha' \bigl( [t'_i,s'_i] \bigr)$ is at most~$L$,
while for $1\le i \le n$ the sets
$\alpha \bigl( [s_{i-1},t_i] \bigr)$
and $\alpha' \bigl( [s'_{i-1},t'_i] \bigr)$ lie in an $L$--neighborhood of
the flat $F_i$.

We will frequently say that paths \emph{$L$--fellow travel relative
to flats} if they $L$--fellow travel relative to some sequence of flats.
\end{definition}

\begin{definition}[Relative Fellow Traveller Property]\label{def:RelFTP}
A space $X$ satisfies the
\emph{Relative Fellow Traveller Property}
if for each choice of constants $\lambda$ and~$\epsilon$
there is a constant $L = L(\lambda,\epsilon,X)$ such that
$(\lambda,\epsilon)$--quasigeodesics in~$X$
with common endpoints $L$--fellow travel relative to flats.
\end{definition}

\section{Quasiflats}
\label{sec:Quasiflats}

This section consists of the proof of Theorem~\ref{thm:Quasiflat}
that quasiflats are close to flats.

\begin{proof}
Let $Q\colon \E^k \to X$ be a $(\lambda, \epsilon)$--quasiflat in~$X$.
Observe that in Euclidean space, any two consecutive sides
of a square form a $\bigl( \sqrt{2},0 \bigr)$--quasigeodesic.
So the image under~$Q$ of any square~$S$ in $\E^k$
can be considered as a pair of $(\lambda',\epsilon')$--quasigeodesics
$\alpha_S$ and $\beta_S$
in~$X$
with common endpoints, where $\lambda'$ and~$\epsilon'$ depend only on
$\lambda$ and~$\epsilon$.

By the Relative Fellow Traveller Property, there is a constant
$L=L(\lambda',\epsilon')$ such that $\alpha_S$ and
$\beta_S$ \  $L$--fellow travel
relative to flats.
But if $S$ is sufficiently large, $\alpha_S$ and $\beta_S$ separate
by more than a distance~$L$ except near their endpoints.
So they must lie within an $L$--neighborhood of some flat~$F_S$.
By isolated flats, this flat lies in a $B$--neighborhood of some
flat $F'_S \in \Set{F}$.

Furthermore, if two sufficiently large squares $S_1$ and~$S_2$
are within a Hausdorff distance~$1$
of each other, then the resulting flats $F'_{S_1}$ and $F'_{S_2}$
must coincide by isolated flats.

To complete the proof, notice that for any given constant~$C$,
one can easily construct a family of
squares $\{S_i\}$ in $\E^k$ each of side length~$C$
with the following two properties.
\begin{enumerate}
\item $\E^k$ is the union of all the squares $S_i$.
\item Any two squares $S$ and~$S'$ can be connected by a finite chain
      \[ S=S_0, S_1, \dots, S_\ell=S' \]
      so that the Hausdorff distance between
      any two consecutive squares $S_i, S_{i+1}$ is at most~$1$.\qedhere
\end{enumerate}
\end{proof}

\section{Equivariant quasi-isometries}
\label{sec:Equivariance}

The goal if this section is to prove Theorem~\ref{thm:EquivariantFTP},
which states that equivariant images
of geodesics lie uniformly close to geodesics.

\begin{definition}[Quasi-equivariance]
Suppose a group $G$ acts by isometries on two metric spaces $X$ and $Y$.
A map $f \colon X \to Y$ is \emph{$\epsilon$--quasi-equivariant}
if for each $g \in G$, the diagram
\[
  \xymatrix{
      X \ar[r]^f \ar[d]^g  &  Y \ar[d]^g\\
      X \ar[r]^f           &  Y
  }
\]
commutes up to a distance~$\epsilon$.  In other words for each $g \in G$
and $x \in X$,
the distance $d \bigl( g(f(x)), f(g(x)) \bigr)$ is less than~$\epsilon$.
\end{definition}

In the sequel, we prove several results about equivariant
quasi-isometries.  In fact, each of these results also holds for
quasi-equivariant quasi-isometries with only trivial modifications to the
proofs (such as introducing extra additive constants).  We will
usually suppress mention of quasi-equivariance
in order to simplify matters slightly.

The proof of Theorem~\ref{thm:EquivariantFTP} uses
the following result, which states that isolated flats and the
Relative Fellow Traveller Property pull back under equivariant
quasi-isometries.  In particular, an equivariant quasi-isometry induces a
one-to-one correspondence between the distinguished families of flats for the
two spaces.

\begin{proposition}\label{prop:Equivariant}
Suppose a group~$G$ acts geometrically on $\CAT(0)$ spaces $X$ and~$Y$.
Suppose further that
\begin{itemize}
   \item $X$ has isolated flats with respect to the family of
      flats $\Set{F}_X$, and
   \item $X$ has the Relative Fellow Traveller Property.
\end{itemize}
Let $\phi \colon Y \to X$ be a $G$--equivariant quasi-isometry.
Then
\begin{enumerate}
\item\label{item:Equivariant:IFP}
   $Y$ has isolated flats with respect to a family of flats
   $\Set{F}_Y$,
\item\label{item:Equivariant:periodic}
   $\phi$ maps the flats of\/ $\Set{F}_Y$ uniformly close to the flats
   of\/ $\Set{F}_X$, inducing a one-to-one correspondence between
   $\Set{F}_Y$ and $\Set{F}_X$, and
\item\label{item:Equivariant:RFTP}
   $Y$ has the Relative Fellow Traveller Property.
\end{enumerate}
\end{proposition}

Curiously, the author does not know of a proof that isolated
flats by itself is a group invariant.  Namely, if a group acts geometrically
on two $\CAT(0)$ spaces and one has isolated flats, must the other
also have isolated flats?  It seems likely that the answer is yes.

\begin{proof}
We first construct the family $\Set{F}_Y$ of flats in~$Y$.
By Theorem~\ref{thm:PeriodicFlats},
each $k$--flat $F \in \Set{F}_X$ is stabilized by a free abelian subgroup
$A \le G$ of rank~$k$.  By the Flat Torus Theorem,
$A$ also stabilizes some $k$--flat $F'$ in~$Y$.
Let $\psi\colon X \to Y$ be a $G$--equivariant quasi-inverse for~$\phi$.
Since $\psi(F)$ and $F'$ are each stabilized by~$A$, it is easy to see that the Hausdorff distance between them is finite.
So $\phi$ sends $\psi(F)$ and $F'$ to a pair of quasiflats which are
each within a finite Hausdorff distance of~$F$.
Applying Theorem~\ref{thm:Quasiflat} to the space~$X$
produces a uniform constant~$B$
which bounds this Hausdorff distance.  Consequently, the Hausdorff distance
between $\psi(F)$ and~$F'$ is
bounded by some other uniform constant $B'$ depending only on the spaces
$X$ and~$Y$ and the constants associated to the maps $\phi$ and~$\psi$.
If we
let the family $\Set{F}_Y$ contain one flat $F' \subseteq Y$ for each flat
$F \in \Set{F}_X$ as constructed above, then
(\ref{item:Equivariant:periodic}) follows immediately.

We now verify that $Y$ has isolated flats with respect to
the family~$\Set{F}_Y$.  Choose a flat $E \subseteq Y$.  By Theorem~\ref{thm:Quasiflat}, the
quasiflat $\phi(E)$ lies within a $D$--neighborhood of some flat $F \in
\Set{F}_X$.  The argument in the previous paragraph shows the existence of a
constant~$D'$ such that $E$ lies in a $D'$--neighborhood of some flat in
$\Set{F}_Y$.
Similarly because $X$ has isolated flats
it is easy to produce a function $\eta \colon \R_{+} \to \R_{+}$
such that for any two flats
$E_1, E_2 \in \Set{F}_Y$ and any constant $C$, the intersection
$\nbd{E_1}{C} \cap \nbd{E_2}{C}$ has diameter bounded by $\eta(C)$.
We have now established~(\ref{item:Equivariant:IFP}).

Finally, let $\alpha$ and~$\alpha'$ be a pair of quasigeodesics in~$Y$
with common endpoints.
Then $\phi\of\alpha$ and $\phi\of\alpha'$ are a pair of
quasigeodesics in~$X$ with common endpoints.
By the Relative Fellow Traveller Property for~$X$,
these quasigeodesics in~$X$ fellow travel relative to
some sequence of maximal flats.
Notice that for any pair of subpaths of $\phi\of\alpha$ and $\phi\of\alpha'$
which are Hausdorff close, the corresponding subpaths of
$\alpha$ and~$\alpha'$ are also Hausdorff close.
On the other hand, given a pair of subpaths $\xi$ and~$\xi'$
which travel far apart but whose endpoints are close, there is some
maximal $k$--flat~$F$ such that $\xi$ and~$\xi'$ both lie close to~$F$.
So in~$Y$ the corresponding subpaths of $\alpha$ and~$\alpha'$ lie close
to the quasiflat $\psi(F)$, which is Hausdorff close to some
$k$--flat~$F'$.  It follows that $\alpha$ and $\alpha'$ fellow travel
relative to some sequence of flats in~$Y$,
establishing~(\ref{item:Equivariant:RFTP}).
\end{proof}

Before proving Theorem~\ref{thm:EquivariantFTP}, we consider the following
special case in which the spaces in question are isometric to Euclidean space.
This case turns out to be quite easy, since the given quasi-isometry
is then close to an affine map.

\begin{lemma}\label{lem:EquivariantZ^n}
Let $\rho_1$ and $\rho_2$ be geometric actions of\/ $\Z^n$ on spaces $F_1$
and~$F_2$ each isometric to Euclidean space~$\E^n$, and let
$\phi \colon F_1 \to F_2$ be a $\Z^n$--equivariant
$(\lambda,\epsilon)$--quasi-isometry.
Choose $\mu$ so that
$F_1$ is contained in a $\mu$--neighborhood of each orbit.
Then the image of a geodesic under~$\phi$ lies within a Hausdorff distance~$L$
of a geodesic, where $L$ depends only on $\lambda$, $\epsilon$, and~$\mu$.
\end{lemma}

\begin{proof}
First choose a basepoint $x\in F_1$, and notice that $\phi$ maps
the orbit of~$x$ to the orbit of $\phi(x)$.
There is an affine map $\psi\colon F_1 \to F_2$ which agrees with $\phi$
on the orbit of~$x$.
But $F_1$ is contained in a $\mu$--neighborhood of this orbit.
So the sup-norm distance between $\psi$ and~$\phi$ is bounded in terms of
$\mu$, $\lambda$, and~$\epsilon$.
But affine maps of Euclidean spaces send lines to lines.
Therefore, $\phi$ sends each geodesic to within a Hausdorff distance~$L$
of a geodesic, where $L=L(\lambda,\epsilon,k)$, as desired.
\end{proof}

\begin{proof}[Proof of Theorem~\ref{thm:EquivariantFTP}]
Proposition~\ref{prop:Equivariant} shows that
$Y$ has the Relative Fellow Traveller Property
and that $\phi$ maps flats in~$X$ uniformly close to
periodic flats in~$Y$.  Pick a quasi-inverse $\psi \colon Y \to X$
for~$\phi$,
and choose a geodesic segment $\alpha$ in~$X$.
Let $\gamma$ be the geodesic in~$Y$ connecting the endpoints of
$\phi\of\alpha$.

It follows from the proof of Proposition~\ref{prop:Equivariant}
that the quasigeodesics $\alpha$ and $\psi\of\gamma$
fellow travel relative to some sequence of flats $(E_1, \dots, E_n)$
in~$\Set{F}_X$, while the quasigeodesics
$\phi\of\alpha$ and $\gamma$ fellow travel
relative to the sequence
$(F_1, \dots, F_n)$ in $\Set{F}_Y$, where $F_i$ is a flat
parallel to the quasiflat $\phi(E_i)$.

Let $\xi$ and $\xi'$ be subsegments of $\alpha$ and $\psi\of\gamma$
which stay
far apart except near their endpoints.  Then the images of $\xi$ and~$\xi'$
lie near some flat $E_i$ in~$X$, while the images of $\phi\of\xi$
and $\phi\of\xi'$ lie near the flat $F_i$.
Recall that $\xi$ is a geodesic parallel to $E_i$, while $\phi\of\xi'$
is a geodesic parallel to $F_i$.

But $\phi$ composed with the orthogonal projection
$\pi_i \colon Y \to F_i$ gives
a quasi-isometry $q\colon E_i \to F_i$,
which is quasi-equivariant with respect to a maximal free abelian subgroup
$A \le G$ that stabilizes both flats.
By Lemma~\ref{lem:EquivariantZ^n}, this quasi-isometry~$q$ maps
geodesics $\epsilon$--close to geodesics
for some constant~$\epsilon$
depending on our choice of flat.  However, without loss of generality
we may assume that
the collections $\Set{F}_X$ and $\Set{F}_Y$ are $G$--equivariant,
so that each consists of only finitely many $G$--orbits of flats.
Therefore, the constants guaranteed by Lemma~\ref{lem:EquivariantZ^n}
are uniformly bounded.
It now follows that $\phi\of\xi$ lies uniformly close
to $\phi\of\xi'$.  Hence, $\phi\of\alpha$ lies close to $\gamma$ as
desired.

The case where $\alpha$ is either a geodesic ray or line follows
easily by a standard argument.
\end{proof}

\section{Limit sets of quasiconvex subgroups}
\label{sec:LimitSets}

In this section, we combine the techniques developed separately
in Sections \ref{sec:QCWellDefined} and~\ref{sec:Equivariance} to prove
Theorem~\ref{thm:LimitSetWellDefined}.

\begin{definition}[Limit Set]
Suppose a group $H$ acts by isometries on a $\CAT(0)$ space~$X$.
Then the \emph{limit set} $\Lambda H$ is the set of accumulation points
in $\boundary X$ of any orbit $Hx$.
Since any two orbits $H x_1$ and $H x_2$ accumulate on the same set,
the choice of basepoint is irrelevant.

If a point~$p$ is the limit of a sequence in $Hx$ that lies within a
finite Hausdorff distance of some geodesic ray, then $p$ is a \emph{conical
limit point} of $H$.  The \emph{conical limit set} $\Lambda_\textrm{c} H$ is
the set of all conical limit points of~$H$.
\end{definition}

Clearly $\Lambda_\textrm{c} H$ lies inside $\Lambda H$.
Furthermore, since $\Lambda H$ and $\Lambda_\textrm{c} H$ are $H$--invariant,
the action of~$H$
on $\boundary X$ restricts to an action of~$H$ on both $\Lambda H$
and $\Lambda_\textrm{c} H$.

The key idea in the proof of Theorem~\ref{thm:LimitSetWellDefined} is
contained in the following lemma, which generalizes the Equivariant Fellow
Traveller Property of Theorem~\ref{thm:EquivariantFTP}.
We will see that the isomorphism $\eta\colon H_1 \to H_2$ induces a map which
sends geodesics in~$X_1$ that lie near the orbit of $H_1$ to geodesics in~$X_2$
that lie near the orbit of~$H_2$.

\begin{lemma}\label{lem:LimitSetEFTP}
Let $\rho_1$ and $\rho_2$ be geometric actions of groups $G_1$ and~$G_2$ on
$\CAT(0)$ spaces $X_1$ and~$X_2$ each with isolated flats
and the Relative Fellow Traveller Property.
Let $H_i \le G_i$ be a subgroup which is $\nu$--quasiconvex with respect
to~$\rho_i$ for some constant~$\nu$.
Let $Y_i = \nbd{H_i x_i}{\nu}$, where $x_i$ is a basepoint in~$X_i$.
Then any isomorphism $\eta \colon H_1 \to H_2$ induces a quasi-isometry
$f \colon Y_1 \to Y_2$ which sends geodesic segments uniformly close to
geodesic segments.
\end{lemma}

\begin{proof}
The map $h \mapsto h(x_i)$ gives a quasi-isometry
$\phi_i \colon H_i \to Y_i$ with quasi-inverse
$\psi_i \colon Y_i \to H_i$.
Since an isomorphism between two finitely generated groups is a quasi-isometry
with respect to any choice of word metrics,
the map $f \colon Y_1 \to Y_2$ given by
$f = \phi_2 \of \eta \of \psi_1$ is a $(\lambda,\epsilon)$--quasi-isometry
for some constants $\lambda$ and~$\epsilon$.
We may also assume that $f$ is $\epsilon$--quasi-equivariant with respect to the
isomorphism~$\eta$;
in other words, for each $h \in H_1$ the diagram
\[
  \xymatrix{
      Y_1 \ar[r]^f \ar[d]^h  &  Y_2 \ar[d]^{\eta(h)}\\
      Y_1 \ar[r]^f           &  Y_2
  }
\]
commutes up to a distance~$\epsilon$.
The map~$f$ has an $\epsilon$--quasi-inverse $g\colon Y_2 \to Y_1$ which is also
$\epsilon$--quasi-equivariant.

We need to see that $f$ maps geodesics uniformly close to
geodesics.  It is enough to prove the result for geodesic segments connecting
two points in the orbit $H_1 x_1$.
So let $\alpha_1 = \bigl[ h(x_1),k(x_1) \bigr]$ for
$h,k \in H_1$.
Then $f \of \alpha_1$ is a $(\lambda,\epsilon)$--quasigeodesic
segment in~$Y_2$.
Perturbing $f\of\alpha_1$ produces a continuous
$(\lambda',\epsilon')$--quasigeodesic $\alpha_2$
with endpoints $\eta(h)(x_2)$
and $\eta(k)(x_2)$ such that the Hausdorff distance
$d_H({f\of\alpha_1}, \alpha_2)$ is bounded in terms of $\lambda$ and~$\epsilon$
and such that the constants $\lambda'$ and~$\epsilon'$ depend only on $\lambda$
and~$\epsilon$.

Since $H_2 x_2$ is a $\nu$--quasiconvex subspace of~$X_2$,
the geodesic segment $\beta_2$ joining the endpoints of~$\alpha_2$ lies
inside~$Y_2$.  As before, $g \of\beta_2$ is within a
bounded Hausdorff distance of a continuous $(\lambda',\epsilon')$--quasigeodesic
$\beta_1$ with the same endpoints as $\alpha_1$, as illustrated in
Figure~\ref{fig:LimitSetEFTP}.

\begin{figure}
\begin{center}
                               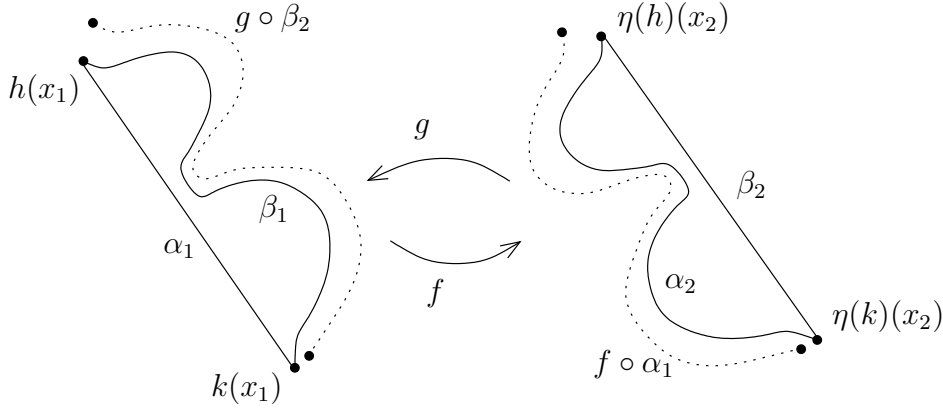
                               \end{center}
\caption{The quasi-isometry~$f$ maps the
pair of quasigeodesics $\alpha_1$ and~$\beta_1$ in~$Y_1$
close to the pair $\alpha_2$ and~$\beta_2$ in~$Y_2$.}
\label{fig:LimitSetEFTP}
\end{figure}

Let $\Set{F}_i$ denote the distinguished family of isolated flats in~$X_i$.
The paths $\alpha_1$ and $\beta_1$ are $(\lambda',\epsilon')$--quasigeodesics
in~$X_1$ with common endpoints.  So
by the Relative Fellow Traveller Property, they $\delta$--fellow travel relative
to the flats in~$\Set{F}_1$ for some $\delta = \delta(X_1,\lambda',\epsilon')$.
Let $\xi$ and~$\xi'$ be subsegments of $\alpha_1$ and $\beta_1$ which stay
far apart except near their endpoints.
Then $\xi$ and~$\xi'$ both lie in a $\delta$--neighborhood of some flat
$F_1 \in \Set{F}_1$.
By Lemma~\ref{lem:PushAcrossOneFlat}, there is a constant $L = L(X_1,\delta)$
such that $\xi$ and~$\xi'$ lie in the $L$--neighborhood of some
subflat $\hat{F}_1 \subseteq F_1$ on which a free abelian subgroup
$A_1 \le H_1$
of rank at least two acts cocompactly by translations.
Furthermore, we may assume that for each $y \in \hat{F}_1$, the orbit
$A_1(y)$ is
$L$--dense in~$\hat{F}_1$.

Since each point of~$\xi$ lies within a uniformly bounded neighborhood
of the orbit $H_1 x_1$ and also within an $L$--neighborhood of $\hat{F}_1$,
there is some $b \in H_1$ so that $b(x_1)$ lies in a uniformly bounded
neighborhood of~$\hat{F}_1$.  So $\hat{F}_1$ lies in a bounded neighborhood
of the orbit $A_1 b(x_1)$.  In particular, there is a quasi-isometry
Furthermore, we may assume that for each $y \in \hat{F}_1$, the orbit
$A_1(y)$ is
$L$--dense in~$\hat{F}_1$.

Since each point of~$\xi$ lies within a uniformly bounded neighborhood
of the orbit $H_1 x_1$ and also within an $L$--neighborhood of $\hat{F}_1$,
there is some $b \in H_1$ so that $b(x_1)$ lies in a uniformly bounded
neighborhood of~$\hat{F}_1$.  So $\hat{F}_1$ lies in a bounded neighborhood
of the orbit $A_1 b(x_1)$.  In particular, there is a quasi-isometry
$\hat{F}_1 \to A_1 b(x_1)$ which moves points by at most a uniformly bounded
amount.
So the composition of quasi-isometric embeddings
\[
  \hat{F}_1 \to A_1 b(x_1) \inclusion Y_1 \to Y_2 \inclusion X_2
\]
is a quasiflat $Q \colon \hat{F}_1 \to X_2$.  Furthermore, if we let
$A_2 = \eta(A_1)$ then the map~$Q$ is quasi-equivariant
with respect to the isomorphism $\eta \big| A_1 \colon A_1 \to A_2$.

By Theorem~\ref{thm:Quasiflat}, there is a universal constant
$D$ so that the quasiflat~$Q$ lies in a $D$--neighborhood
of some flat $F_2 \in \Set{F}_2$ which is stabilized by~$A_2$.
Since $Q$ is quasi-equivariant and $A_2$ acts on~$F_2$ by translations,
$Q$ must lie inside a uniformly bounded neighborhood of some subflat
$\hat{F}_2 \subseteq F_2$ on which $A_2$ acts cocompactly.

As in the proof of Theorem~\ref{thm:EquivariantFTP}, projecting $Q$
onto~$\hat{F}_2$ gives an equivariant quasi-isometry $\hat{F}_1 \to \hat{F}_2$.
Such a map sends geodesics uniformly close to geodesics by
Lemma~\ref{lem:EquivariantZ^n}.  Since $\xi$ is a geodesic segment in the
$L$--neighborhood of~$\hat{F}_1$, its image $f \of \xi$ in~$X_2$ lies close
to a geodesic.  But the endpoints of $f \of \xi$ are close to the
geodesic~$\beta_2$.
It now follows that $f$ maps the entire geodesic~$\alpha_1$
into a uniformly bounded neighborhood of~$\beta_2$.
The uniform bound in question depends only on our original choice of
quasi-isometric embeddings $H_i \to X_i$ and on the given isomorphism
$\eta\colon H_1 \to H_2$.
\end{proof}

\begin{proof}[Proof of Theorem~\ref{thm:LimitSetWellDefined}]
Fix a basepoint $X_i \in X_i$, and consider the quasi-isometric embedding
$H_i \to X_i$ given by $h \mapsto h(x_i)$.
By Lemma~\ref{lem:LimitSetEFTP}, any isomorphism $\eta\colon H_1 \to H_2$
induces a quasi-isometry~$f$ from
\[
   Y_1 = \bignbd{H_1(x_1)}{2\nu} \quad \text{to} \quad
   Y_2 = \bignbd{H_2(x_2)}{2\nu}
\]
which sends geodesic segments uniformly close to
geodesics.  It follows easily that $f$ maps geodesic rays uniformly close
to geodesic rays.
Thus we have a one-to-one correspondence between rays in~$Y_1$ and rays
in~$Y_2$.  To complete the proof we need to see that every point
of $\lambda H_i$ can be represented by a ray in~$Y_i$.

Consider a sequence $\bigl\{ h_j(x_i) \bigr\}$ limiting to a point of
$\boundary X_i$ as $j \to \infty$.
Extract a subsequence so that the segments
$\bigl[ h_1(x_i),h_j(x_i) \bigr]$
converge pointwise to a geodesic ray~$c$ based at $h_1(x_i)$.
By quasiconvexity, each segment
$\bigl[ h_1(x_i),h_j(x_i) \bigr]$ lies inside
the $\nu$--neighborhood of the orbit $H_i(x_i)$.
So the limiting ray~$c$ lies inside $Y_i$.
Therefore, every point of $\Lambda H_i$ is represented by a geodesic ray
inside~$Y_i$ based at the point $h_1(x_i)$.
\end{proof}

\section{Geometrically finite groups}
\label{sec:GeomFinite}

This section is devoted to proving Theorem~\ref{thm:GeomFiniteRFTP}.
We begin by considering the finite volume case.

A \emph{truncated hyperbolic space} is a subspace of $\Hyp^n$ obtained
by removing a collection of disjoint open horoballs
and endowing the resulting subset with the induced length metric.
Every truncated hyperbolic space is a complete $\CAT(0)$ space
(\cite[II.11.27]{BH99}).
A discrete subgroup
$\Gamma \le \Isom(\Hyp^n)$ with finite covolume acts cocompactly
on a truncated space obtained by removing a $\Gamma$--equivariant family
of horoballs centered at the parabolic fixed points of~$\Gamma$
(\cite{GarlandRaghunathan70}, see also \cite[\S 4.5]{Thurston97}).

\begin{proposition}\label{prop:Truncated}
Let $X \subset \Hyp^n$ be any truncated hyperbolic space.
Then $X$ has isolated flats.
\end{proposition}

\begin{proof}
We may assume $n \ge 3$, since otherwise $X$ is $\delta$--hyperbolic.
For every deleted horoball, the bounding horosphere is
isometric to $\E^{n-1}$.  Since $X$ is locally isometric to $\Hyp^n$
away from these flats, these horospheres are the only flats in~$X$.

To verify that $X$ has isolated flats, we need to bound the diameter
$D(k)$ of the intersection of $k$--neighborhoods of any two distinct flats.
Notice that a tubular neighborhood of a horoball is again a horoball.
Furthermore the diameter of the intersection decreases monotonicaly
as a function of the distance between the two flats.
So it suffices to consider the case where the two horoballs are tangent
at a single point, in which case it is clear that the diameter obtained
is finite and depends only on~$k$.
\end{proof}

A geometrically finite group $\Gamma\le \Isom(\Hyp^n)$
acts geometrically on a \emph{truncated convex hull} obtained as follows.
Let $\Lambda$ be the limit set of $\Gamma$ in $\boundary\Hyp^n$,
and let $\Hull(\Lambda) \subseteq \Hyp^n$ be the hyperbolic convex hull
of $\Lambda$.  If $\Gamma$ is geometrically finite, then
there is a $\Gamma$--equivariant collection of disjoint open horoballs,
with union~$U$, centered at the parabolic fixed points of~$\Gamma$, such that the action of $\Gamma$ on the truncated convex hull
$Y = \Hull(\Lambda) \cap (\Hyp^n \setminus U)$
is properly discontinuous and cocompact (\cite{Bowditch93}).

If the horoballs in~$U$ are chosen sufficiently small, then the truncated
convex hull is a convex subspace of the truncated hyperbolic space,
and hence is $\CAT(0)$.  This fact seems to be well-known,
though the only explicit reference the author has found in the literature
is Exercise II.11.37(2) of \cite{BH99}.
We direct the reader towards
the lemma in \cite[\S 1.7]{CullerShalen92}, which is a key step
in proving this exercise.
This lemma is proved by Culler--Shalen only in the
three-dimensional setting, but the generalization to higher dimensions
is straightforward.  Henceforth, we assume that all truncated convex
hulls are $\CAT(0)$.

\begin{proposition}[Geometrically finite hyperbolic]
\label{thm:GeometricallyFinite}
Let $\Gamma$ be any geometrically finite subgroup of\/ $\Isom(\Hyp^n)$.
Then the associated truncated convex hull~$Y$ has isolated flats.
\end{proposition}

\begin{proof}
As in the proof of Proposition~\ref{prop:Truncated},
the only flats of~$Y$ are contained in the bounding horospheres.
The stabilizer in~$\Gamma$ of each horosphere~$S$
is virtually abelian of rank $k$ for some $k<n$.
The intersection of~$S$ with $\Hull(\Lambda)$
is isometric to a product $F \times Z$ with $F$ isometric to~$\E^k$
and $Z$ a compact convex subset of $\E^{n-k-1}$.
Let $z\in Z$ be the circumcenter of~$Z$ (see \cite[II.2.7]{BH99}).
Define $\Set{F}$ to be the set of all flats
$F \times \{z\}$ whose stablilizer in~$\Gamma$ has rank at least two.
By construction, $\Set{F}$ is invariant under $\Isom(Y)$.
Since the horospheres~$S$ are isolated, it follows that
$\Set{F}$ is isolated as well.
Furthermore, since $\Gamma$ has only finitely many conjugacy classes
of maximal parabolic subgroups, it is easy to see that
each flat in~$Y$ lies in a universally bounded neighborhood
of some flat in~$\Set{F}$.
\end{proof}

In order to prove Theorem~\ref{thm:GeomFiniteRFTP},
all that remains is to establish the Relative Fellow Traveller Property
for the truncated convex hull, which follows easily from
the following result due to Epstein
(\cite[Theorem~11.3.1]{ECHLPT92}).

\begin{theorem}[Quasigeodesics outside horoballs]
\label{thm:EpsteinHoroballs}
Let $\lambda \ge 1$ and $\epsilon \ge 0$ be fixed real constants.
Then there is a positive real number~$\ell$, depending only on
$k$ and~$\epsilon$, with the following property.
Let $r > 3\ell$.
Let $U$ be a union of disjoint horoballs in $\Hyp^n$, such that
any two components of~$U$ are a distance at least~$r$ apart,
and let $X$ be the truncated space $\Hyp^n \setminus U$.
Let $\alpha\colon [a,b] \to X$ be a
$(\lambda,\epsilon)$--quasigeodesic in~$X$.
Let $\phi$ be the hyperbolic geodesic from $\alpha(a)$ to $\alpha(b)$.
Then the union of the $\ell$--neighborhood
of~$U$ and the $\ell$--neighborhood of~$\phi$ contains the image of~$\alpha$.
\end{theorem}

\begin{proof}[Proof of Theorem~\ref{thm:GeomFiniteRFTP}]
Consider the truncated hyperbolic space $X = \Hyp^n \setminus U$
and truncated convex hull $Y = \Hull(\Lambda) \cap X$
associated to~$\Gamma$.
Shrink the horoballs in~$U$ equivariantly so that they
satisfy the hypothesis of Theorem~\ref{thm:EpsteinHoroballs}.

Since $Y$ is convex in~$X$, quasigeodesics in~$Y$
are also quasigeodesics in~$X$.
So the result follows from the fact that,
for each bounding horosphere~$S$, the intersection $S \cap Y$
lies uniformly close to either a flat or a geodesic line.
\end{proof}

\bibliographystyle{alpha}
\bibliography{chruska}

\end{document}